\documentclass[12pt]{amsart}
\usepackage{amssymb}
\usepackage[all]{xy}

\usepackage{a4wide}

\newtheorem{thm}{Theorem}[section]
\newtheorem{lem}[thm]{Lemma}

\newtheorem{prop}[thm]{Proposition}

\newtheorem{cona}[thm]{The generalized Auslander conjecture}
\newtheorem*{prob*}{Open problem}

\theoremstyle{definition}

\newtheorem{defi}[thm]{Definition}

\theoremstyle{remark}

\newtheorem*{rem*}{Remark}

\DeclareMathOperator{\id}{id}
\DeclareMathOperator{\s}{span}

\DeclareMathOperator{\Aff}{Aff}

\newcommand{\kringel}{\mathbin{\raise1pt\hbox{$\scriptstyle\circ$}}}
\newcommand{\pkt}{\mathbin{\raise0pt\hbox{$\scriptstyle\bullet$}}}

\newcommand{\C}{\mathbb{C}}

\newcommand{\N}{\mathbb{N}}

\newcommand{\R}{\mathbb{R}}

\newcommand{\Z}{\mathbb{Z}}

\newcommand{\Der}{\mathop{\rm Der}}
\newcommand{\Aut}{\mathop{\rm Aut}}
\newcommand{\GL}{\mathop{\rm GL}}
\newcommand{\SL}{\mathop{\rm SL}}

\newcommand{\La}{\mathfrak{a}}
\newcommand{\Lb}{\mathfrak{b}}

\newcommand{\Ld}{\mathfrak{d}}
\newcommand{\Le}{\mathfrak{e}}

\newcommand{\Lg}{\mathfrak{g}}
\newcommand{\Lh}{\mathfrak{h}}

\newcommand{\Ln}{\mathfrak{n}}

\newcommand{\abs}[1]{\lvert#1\rvert}

\newcommand{\norm}[1]{\lVert #1 \rVert}

\newcommand{\al}{\alpha}
\newcommand{\be}{\beta}
\newcommand{\ga}{\gamma}
\newcommand{\de}{\delta}
\newcommand{\ep}{\varepsilon}

\newcommand{\la}{\lambda}

\newcommand{\ov}{\overline}

\newcommand{\ra}{\rightarrow}

\renewcommand{\phi}{\varphi}

\begin{document}


\title[The Auslander conjecture]{The Auslander conjecture for NIL-affine crystallographic groups}

\author[D. Burde]{Dietrich Burde}
\author[K. Dekimpe]{Karel Dekimpe}
\author[S. Deschamps]{Sandra Deschamps}
\address{Katholieke Universiteit Leuven\\
Campus Kortrijk\\
8500 Kortrijk\\
Belgium}
\date{\today}
\email{dietrich.burde@univie.ac.at}
\email{karel.dekimpe@kulak.ac.be, sandra.deschamps@kulak.ac.be}

\begin{abstract}
We study subgroups $\Gamma$ in $\Aff (N)=N\rtimes \Aut (N)$ acting
properly discontinuously and cocompactly on $N$. Here $N$ is a simply connected, connected real nilpotent
Lie group of finite dimension $n$. This situation is a natural generalization of the so-called affine
crystallographic groups. We prove that for all dimensions $1\le n\le 5$
the generalized Auslander conjecture holds, i.e., that such subgroups are virtually
polycyclic.
\end{abstract}

\maketitle

\section{Introduction}
A classical crystallographic group is a discrete subgroup of ${\rm
Isom}(\R^n)$. Such groups act properly discontinuously and
cocompactly on $\R^n$. The structure of such groups is well known
by the three Bieberbach theorems (\cite{CHAR}, \cite{WOLF}). In
fact, all these groups are finitely generated virtually abelian.

As a generalization of this concept, one also studies affine
crystallographic groups. These are subgroups of $\Aff(\R^n)=\R^n\rtimes \GL_n(\R)$
acting crystallographically (by which we will always mean properly
discontinuously and cocompactly) on $\R^n$. The structure of these
affine crystallographic groups is not at all as well known as in
the case of the classical crystallographic groups and two main
questions have played a major role in the study of these groups:
\begin{itemize}
\item In 1977 Milnor (\cite{MIL}) asked whether any torsionfree polycyclic-by-finite group could be
realized as an affine crystallographic group. And conversely
\item In 1964 Auslander conjectured (\cite{AUS}) that any affine crystallographic group is
virtually solvable (and hence polycyclic-by-finite).
\end{itemize}
It follows that a positive answer to both questions implies a
complete understanding of the affine crystallographic groups.
However, it has been shown that not all torsionfree
polycyclic-by-finite groups do occur as an affine crystallographic
group, thus answering negatively Milnor's question (\cite{BENO},
\cite{BG}, \cite{BU}).

On the other hand, Auslander's conjecture is still open and is
only known to be true in dimensions $n\le 6$ (\cite{AMS2}, see
also \cite{ABELS} for a survey on the Auslander conjecture).

The negative answer to Milnor's question, suggests to consider
more general classes of crystallographic groups (e.g.\ polynomial
crystallographic groups as in \cite{DEI}). In this paper we are
considering crystallographic subgroups in $\Aff(N)$ for a simply
connected, connected nilpotent Lie group $N$. Here $N$ is
diffeomorphic to some $\R^n$ and $\Aff (N)=N\rtimes \Aut (N)$ acts
on $N$ via
\begin{align}\label{action1}
(n,\al)\cdot m & = n \al(m) \quad \forall \; (n,\al)\in \Aff(N),\;
m\in N.
\end{align}
A crystallographic subgroup of $\Aff(N)$ is a subgroup acting crystallographically
on $N$. We will call such a group a NIL-affine crystallographic group.\\
The notation $\Aff(N)$ makes sense, since this group is really the
group of connection preserving diffeomorphisms of $N$ for any left
invariant affine connection on $N$ (\cite{KT}). In this sense,
studying the NIL-affine crystallographic groups is really a very
natural generalization of the affine crystallographic groups.
Moreover, the analogue of Milnor's question does hold in this
case:
\begin{thm}$($\cite{DE1}$)$
Let $\Gamma$ be a torsionfree polycyclic-by-finite group. Then there exists a connected and simply connected
nilpotent Lie group $N$ and an embedding $\rho:\Gamma\rightarrow \Aff(N)$, such that
$\rho(\Gamma)$ is a crystallographic subgroup of $\Aff(N)$.
\end{thm}
Conversely, just as in the case of affine crystallographic groups,
it is now very natural to ask whether every NIL-affine
crystallographic group is virtually solvable:
\begin{cona}
Let $N$ be a connected and simply connected nilpotent Lie group and let
$\Gamma \subseteq \Aff(N)$ be a group acting crystallographically on $N$.
Then $\Gamma$ is virtually polycyclic.
\end{cona}
We expect the answer to be positive since the question is closely
related to the original Auslander problem. Note that a positive
answer to this generalized Auslander problem would imply a
complete algebraic description of the class of NIL-affine
crystallographic groups, or stated otherwise, would provide a
complete geometric description of the class of polycyclic-by-finite groups.\\
We will prove this conjecture for all $N$ with $\dim N\le 5$.
Moreover, in the next section we will show how several concepts of
classical affine geometry on $\R^n$ can be translated to the
nilpotent case, giving more indications of a close relation
between the generalized Auslander conjecture and the original one
and thus providing even more evidence for a positive answer.

\section{Subgroups of $\Aff(N)$ not acting properly
discontinuously}\label{crit}

In this section we generalize the criterium of \cite{AMS} for the
failure of proper discontinuity for certain subgroups in $\Aff(\R^n)$
to the case of subgroups in $\Aff(N)$.
This means that we follow on one hand very closely the construction given
in \cite{AMS}, but on the other hand really have to develop some basics
for affine geometry in a nilpotent Lie group $N$. Here we encounter new
problems, because $N$ will not be commutative in general.

\medskip

In establishing this criterium, we will exploit the structure of the
{\em linear} parts of the affine motions involved. Let us make this
more precise in what follows.\\
Given $\al\in \Aut(N)$, let $\al_{\ast}$ denote the induced automorphism
of the Lie algebra $\Ln$ of $N$. We have the identification
$\R^n\equiv \Ln \xrightarrow{\exp} N$, $\exp\kringel\al_{\ast}=\al\kringel \exp$
and $\Aut(N)\cong \Aut(\Ln)$.
Denote by
\[
\ell\colon N\rtimes \Aut (N) \ra \Aut(\Ln):\quad (n,\al)\mapsto \al_{\ast}
\]
the projection of $g=(n,\al)$ to its linear part
$\ell(g)=\al_{\ast}$. Then we have for any $g\in \Aff(N)$ the
vector space decomposition of $\Ln$ into a direct sum of
$\ell(g)$-invariant subspaces

\begin{align}\label{vsd1}
\Ln & = \Ln^{-}(g)\oplus \Ln^{0}(g) \oplus \Ln^{+}(g)
\end{align}

where the spaces $ \Ln^{-}(g),\, \Ln^{0}(g),\, \Ln^{+}(g)$ are determined
by the following conditions. Their sum is $\Ln$ and
all eigenvalues $\la$ of the restriction $\ell(g)_{\mid \Ln^{-}(g)}$ satisfy
$\abs{\la}<1$, all eigenvalues of $\ell(g)_{\mid \Ln^{0}(g)}$ satisfy $\abs{\la}=1$
and all eigenvalues of $\ell(g)_{\mid \Ln^{+}(g)}$ satisfy $\abs{\la}>1$.
The decomposition is not only a vector space decomposition but also a decomposition
as Lie algebras. In fact, we will need that the subspaces  $ \Ln^{-}(g),\, \Ln^{0}(g),\,
\Ln^{+}(g)$ and the two direct sums

\begin{align}\label{vsd2}
\Ld^{-}(g) & =  \Ln^{-}(g)\oplus \Ln^{0}(g)\\
\Ld^{+}(g) & = \Ln^{+}(g)\oplus \Ln^{0}(g)
\end{align}

are Lie subalgebras of $\Ln$. However,  we will prove this in 
Lemma~\ref{lemma2.1}. So we can also fix notations for
the corresponding Lie subgroups of $N$:
\begin{align*}
N^-(g) & = \exp(\Ln^-(g))\\
N^0(g) & = \exp(\Ln^0(g))\\
N^+(g) & = \exp(\Ln^+(g))\\
D^-(g) & = \exp(\Ld^-(g))\\
D^+(g) & = \exp(\Ld^+(g)).
\end{align*}

\begin{lem}\label{lemma2.1}
Let $\Lg$ be a $n$-dimensional real Lie algebra and $\al\in \Aut(\Lg)$.
Then there exists a basis $(v_1,\ldots,v_n)$ of $\Lg$ such that $\al$
has the following block form with respect to this basis
\[
\al=\begin{pmatrix}
A & \vrule & 0 & \vrule & 0 \\
\hline \\[-0.48cm]
0 & \vrule & B & \vrule & 0 \\
\hline \\[-0.48cm]
0 & \vrule & 0 & \vrule & C
\end{pmatrix}
\]
where $A\in M_k(\R)$ such that all eigenvalues $\la$ satisfy $\abs{\la}<1$,
$B\in M_l(\R)$ such that all eigenvalues $\la$ satisfy $\abs{\la}=1$, and
$C\in M_m(\R)$ such that all eigenvalues $\la$ satisfy $\abs{\la}>1$.
If
\begin{align*}
\Lg^{-}(\al) & =  \s \{ v_1,\ldots , v_k\} \\
\Lg^{0}(\al) & =  \s \{ v_{k+1},\ldots , v_{k+l}\} \\
\Lg^{+}(\al) & =  \s \{ v_{k+l+1},\ldots , v_{k+l+m} \}
\end{align*}
with $k+l+m=n$ then the subspaces
$\Lg^{-}(\al),\, \Lg^{0}(\al),\, \Lg^{+}(\al)$ and $\Lg^{-}(\al)\oplus \Lg^{0}(\al),
\, \Lg^{+}(\al)\oplus \Lg^{0}(\al)$ are Lie subalgebras of $\Lg$.\\
Moreover, $\Lg^+(\al)$ (resp.\ $\Lg^-(\al)$) is an ideal of
$\Lg^+(\al)\oplus\Lg^{0}(\al)$ (resp.\
$\Lg^-(\al)\oplus\Lg^{0}(\al)$).
\end{lem}

\begin{proof}
The first assertion follows easily by using the real canonical Jordan
form for $\al$. Now let $\al=\al_s \al_u$ be the multiplicative Jordan
decomposition of $\al$. Here $\al_s$ is a semisimple automorphism,
$\al_u$ is a unipotent automorphism and $\al_s \al_u=\al_u \al_s$. We have
$\al_s, \al_u\in \Aut(\Lg)$ since $\Aut(\Lg)$ is a linear algebraic group.
In fact, $\al_s$ and $\al_u$ are represented by block matrices as above,
\[
\al_s=\begin{pmatrix}
A_s & 0   & 0  \\
0   & B_s & 0  \\
0   & 0   & C_s
\end{pmatrix},\quad
\al_u=\begin{pmatrix}
A_u & 0   & 0  \\
0   & B_u & 0  \\
0   & 0   & C_u
\end{pmatrix}
\]
where the subscript $s$ means that we take the semisimple part of
the matrix, and the subscript $u$ stands for the unipotent part.
Note that $\Lg^{\ep}(\al)=\Lg^{\ep}(\al_s)$ for $\ep=-,0,+$. We
may assume that the matrices appearing in the representation of
$\al_s$ have diagonal form. Otherwise we may pass to the
complexification of $\Lg$ where we can diagonalize. Now a direct
calculation finishes the proof. Let us first check that the space
$\Lh=\Lg^{-}(\al)\oplus \Lg^{0}(\al)$ is a Lie subalgebra of
$\Lg$. All other cases are analogous. $\Lh$ is spanned by all
eigenvectors corresponding to an eigenvalue $\lambda$ with
$|\lambda|\leq 1$. So, suppose that $\{ v_1,\ldots v_{k+l}\}$ is a
basis of $\Lh$ such that $\al_s(v_i)=\la_i v_i$ for all $i$
$(|\lambda_i|\leq1)$. For any $v_i,v_j \in \Lh$, we compute that
\begin{align*}
\al_s([v_i,v_j]) & = [\al_sv_i,\al_sv_j] = [\la_iv_i,\la_jv_j] =\la_i\la_j [v_i,v_j].
\end{align*}
Now, clearly $\abs{\la_i\la_j}\le 1$,
from which it follows that $[v_i,v_j]\in \Lh$, proving that $\Lh$ is a Lie
subalgebra of $\Lg$. \\
Checking that $\Lg^{-}(\al)$ is an ideal of $\Lh$ can be done
similarly, which finishes the proof.
\end{proof}

To be able to generalize the ideas of \cite{AMS}, we must be able
to talk about {\em affine subspaces} of $N$. Therefore, we define
a line in the Lie group $N$ as a left coset of a $1$-parameter
subgroup. In other words, $L=m\cdot \exp (tA)$ for some $A\in \Ln$
with $A\neq 0$. We say that this line is parallel to the Lie
subalgebra $\s \{A\}$. More generally, we have the following
definition.

\begin{defi}
For a Lie subalgebra $\Lh\subseteq \Ln$ we define an affine
subspace of $N$ to be any left coset of the form $m\cdot
\exp(\Lh)$ ($m\in N$). We say that this affine subspace is {\it
parallel} to $\Lh$.
\end{defi}

Define the following two subsets of $\Aff(N)$:
\begin{align*}
\Omega & =\{ g\in \Aff(N) \mid \dim \Ln^0 (g)=1,\; \ell(g)_{\mid \Ln^{0}(g)} =\id \} \\
\Omega_0 & = \{ g\in \Omega \mid g\cdot n \neq n \quad \forall \, n\in N \}.
\end{align*}
The elements of  $\Omega$ are called  pseudohyperbolic, and
$\Omega_0$ consists of fixed-point-free elements of $\Aff(N)$ inside $\Omega$.
We will study the action of a group generated by pseudohyperbolic elements.
The following lemma is needed to be able to describe the behaviour of
the action of a pseudohyperbolic element $g$ and its iterates $g^n$.

\begin{lem}\label{lemma2.3}
Suppose that $\al\in \Aut(N)$ satisfies the following condition:
if $1$ is an eigenvalue of $\al_{\ast}\in \Aut(\Ln)$ then its
geometric and algebraic multiplicity coincide. Let
$\Le=Eig(\al_{\ast},1)$ be the eigenspace to the eigenvalue $1$,
and $E=\exp(\Le)$. For fixed $n\in N$ define a map $\phi\colon
N\ra N$ by
\[
m\mapsto m^{-1}n\al(m).
\]
Then there exists an $m\in N$ such that $\phi(m)\in E$. This $m$ is uniquely
determined modulo $E$.
\end{lem}

\begin{proof}
We prove the result by induction on the nilpotency class $c$ of $N$.
If $N$ is abelian, i.e., for c=1 and $N=\R^n$, we have
\begin{align}\label{abelian}
\phi(m)& =-m+n+\al(m)=(\al-\id)\cdot m +n.
\end{align}
Let $(v_1,\ldots ,v_n)$ be a basis of $\R^n$ such that the first
$k$ vectors form a basis of $E$. Then $\al$  is represented by the
matrix
\[
\al=\begin{pmatrix}
I & \vrule & \ast \\
\hline \\[-0.48cm]
0 & \vrule & A
\end{pmatrix}
\]
where $A\in M_{n-k}(\R)$ has no eigenvalues equal to $1$. Hence the block matrices
on the diagonal of $\al-\id$ are $0$ and $A-I$, the latter being invertible.
Hence the matrix equation $(A-I)x+b=0$ for any $b\in \s \{v_{k+1},\ldots ,v_n\}$ has a
unique solution. This means that there is an $m\in N$
such that the last $n-k$ components of the vector $\phi(m)=(\al-\id)\cdot m +n$
are zero, i.e., such that $\phi(m)\in E$.
Moreover the last $n-k$ components of $m$ are uniquely determined.
Hence if $m'\in N$ is another element satisfying $\phi(m')\in E$, then $m=m'+e$
with some $e\in E$. \\
Now suppose that $c>1$ and that the lemma is true for lower nilpotency classes.
Let $Z$ be the center of $N$ and define $\ov{\phi}\colon N/Z \ra N/Z$ by
\[
\ov{m}\mapsto \ov{m^{-1}n\al(m)}=(\ov{m})^{-1}\ov{n}\ov{\al(m)}
\]
where $\ov{\al}\colon  N/Z \ra N/Z$ given by $\ov{m}\mapsto \ov{\al(m)}$
is an automorphism of $N/Z$. (Here we use the bar to denote the natural
projection $N\ra N/Z$). Note that $\ov{\phi}$ is well-defined and that
$\ov{\al}$ satisfies the assumption of the lemma on the eigenvalue $1$. Hence
we can apply the induction hypothesis, and there is an $\ov{m}\in N/Z$ such that
$(\ov{m})^{-1}\ov{n}\ov{\al(m)}\in EZ/Z \subseteq N/Z$. Thus (for any given
lift $m$ of $\ov{m}$) we may write
\[
m^{-1}n\al(m)=e_1z_1
\]
with some $e_1\in E$ and $z_1\in Z$. It follows that for any $z\in
Z$
\begin{align*}
\phi(mz) & = z^{-1}m^{-1}\cdot n\cdot \al(m)\al(z)\\
  & = m^{-1}n\al(m)\cdot z^{-1}\al(z)\\
  & = e_1z_1\cdot  z^{-1}\al(z)\\
  & = e_1\cdot z^{-1}z_1\al(z).
\end{align*}
Since $\al_{\mid Z}\in \Aut (Z)$ satisfies the eigenvalue $1$
condition of the lemma we can again apply the induction
hypothesis, with $N=Z$ and $n=z_1$. Hence we find a $z_0\in Z$
such that $z_0^{-1}z_1\al(z_0)\in E\cap Z$. It follows that
\[
\phi(mz_0)=e_1\cdot z_0^{-1}z_1\al(z_0)\in E.
\]
The uniqueness up to $E$ also follows by induction on the nilpotency class of $N$.
\end{proof}

Now, we apply this lemma to obtain some information about the action of a
pseudohyperbolic element.

\begin{prop}\label{1line}
For any $g\in \Omega$ there is exactly one $g$-invariant line $C_g$ parallel to
$\Ln^0 (g)$.
\end{prop}

\begin{proof}
A line being parallel to $\Ln^0 (g)$ is of the form $m\cdot \exp(tA)$
with some $m\in N$ and an $A$ satisfying $\Ln^0 (g)=\s \{A\}$.
It is $g$-invariant if and only if there is a function $s\colon \R \ra \R$ such that
\begin{equation}\label{LL}
g(m\cdot \exp(tA))=m\cdot \exp(s(t)\cdot A).
\end{equation}
Writing $g=(n,\al)$ and applying \eqref{action1} we obtain
\begin{align}
g(m\cdot \exp(tA))& =  n\cdot \al(m)\cdot \al(\exp(tA)) \nonumber \\
 & =  n\cdot \al(m)\cdot \exp(tA). \label{RL}
\end{align}
Now, equating (\ref{LL}) to (\ref{RL}) leads to
$m^{-1}n\cdot\al(m)=\exp((s(t)-t)A)$,
which is an element of the $1$-dimensional Lie group $N^0(g)=\exp(\Ln^0(g))$.
We conclude that there exists a $g$-invariant line parallel to
$\Ln^0 (g)$ if and only if there is an $m\in N$ such that
$m^{-1}n\cdot\al(m)\in N^0(g)$ (and then $s(t)-t$ has to be constant).
But this follows by Lemma $\ref{lemma2.3}$. We have $\Le=\Ln^0(g)$ and
$E=N^0(g)$. Moreover $\al$ satisfies the eigenvalue $1$ condition. Hence there is an
$m\in N$ and some $c\in \R$ with $m^{-1}n\cdot\al(m)=\exp(cA)\in N^0(g)$.
The line $m\cdot \exp(tA)$ is $g$-invariant and we have
\begin{equation}\label{translation}
g(m\cdot \exp(tA))=m\cdot \exp((c+t)A)=m\cdot \exp(tA)\exp(cA).
\end{equation}
Since $m$ was unique up to $N^0(g)$, another choice of $m$ yields
the same line. Indeed, if $m'=mn_1$ with some $n_1=\exp(c_1A)\in
N^0(g)$ then
\[
m'\exp(tA)=mn_1\exp (tA)=m\exp((c_1+t)A).
\]
Hence this line is unique.
\end{proof}

The proposition above does not only give us a $g$-invariant line
$C_g$, for any $g\in \Omega$, but equation $(\ref{translation})$ shows that
the action on any point $x$ of this line is by means of a constant translation.
We define the {\it translational part} $\tau(g)$ of $g$ by $gx=x\tau(g)$
(where $x$ is any point $x\in C_g$).
It holds that $\tau(g)\neq 1$ if and only if $g\in \Omega_0$. For $m\in \Z, m\neq 0$ we have
\begin{align*}
C_{g^m} & = C_g\\
\tau(g^m) & = \tau(g)^m.
\end{align*}
Let $T(g)=\log (\tau(g))\in \Ln^0(g)$. If $g\in \Omega_0$, then $T(g)\neq 0$,
and hence every $x\in \Ld^+(g)$
has a unique decomposition
\[
x=\la(x)T(g)+a(x)
\]
where $\la(x)\in \R$ and $a(x)\in \Ln^+(g)$.
We call $x\in \Ld^+(g)$ {\it positive with respect to} $g\in \Omega_0$, if $\la(x)>0$.
We will write $x\succ _g 0$.

\begin{defi}
Two elements $g_1,g_2\in \Omega_0$ will be called {\it transversal}, if
\[
\Ln=\Ln^+(g_1)\oplus \Ld^+(g_2)=\Ld^+(g_1)\oplus \Ln^+(g_2).
\]
\end{defi}

It is easy to see that $g_1,g_2\in \Omega_0$ are transversal if and only if
\begin{align*}
\Ln^+(g_1)\oplus \Ln^+(g_2) \oplus (\Ld^+(g_1) \cap \Ld^+(g_2)) & = \Ln\\
\mbox{and }\dim (\Ld^+(g_1) \cap \Ld^+(g_2)) & = 1.
\end{align*}

Let
\begin{align*}
S_{g_1} & = \{ x\in \Ld^+(g_1) \cap \Ld^+(g_2) \mid x\succ_{g_1}0 \}\\
S_{g_2} & = \{ x\in \Ld^+(g_1) \cap \Ld^+(g_2) \mid x\succ_{g_2}0
\}.
\end{align*}

\begin{defi}
For two transversal elements $g_1,g_2\in \Omega_0$ we say that
they form a {\it positive pair} if $S_{g_1}=S_{g_2}$.
\end{defi}

For $g\in \Omega$ note that $C_g=m\cdot N^0(g)$.

For any $x\in N$ let \[B_g^+(x) = x\cdot N^+(g).\]

We will also use
\begin{align*}
E_g^+ & = m\cdot D^+(g) = C_g\cdot  D^+(g)\\
E_g^- & = m\cdot D^-(g) = C_g\cdot  D^-(g).
\end{align*}

Note that $\Ln^+(g)$ is a Lie ideal in $\Ld^+(g)$ by Lemma $\ref{lemma2.1}$.
Hence $N^+(g)$ is a normal subgroup in $D^+(g)$. Since the intersection
of $N^+(g)$ and $N^0(g)$ is trivial, every element $x\in D^+(g)$ can be written
as
\[
x=n^0n^+
\]
with unique elements $n^0\in N^0(g)$ and $n^+\in N^+(g)$.
Let $x\in E_g^+$. Then it is easy to see that $B_g^+(x)\cap C_g$ consists of
exactly $1$ point: if $x=m\cdot n^0n^+$ and $C_g=m\cdot N^0(g)$ then
\[
m\cdot n^0=x\cdot (n^+)^{-1}
\]
is this point. Thus we can define a projection
\[
P_g\colon E_g^+ \ra C_g
\]
by the equality $B_g^+(x)\cap C_g=\{P_g(x)\}$ for $x\in E_g^+$.
The subgroup $N^+(g)$ is $\al$-invariant, where $g=(n,\al)$.
Therefore, if $x=m\cdot n^0n^+ \in E_g^+$ then
\begin{equation}
P_g(gx) = m\cdot n^0\tau(g)=P_g(x)\tau(g).
\end{equation}

We are now ready to prove the obstruction criterium to proper discontinuity.

\begin{prop}\label{prop2.7}
Assume that $g_1,g_2\in \Omega_0$ form a positive pair. Then there exists
a compact set $K\subset N$ and two sequences $\{s_i\}, \{t_i\}$ of positive
integers such that

\[
\lim_{i\to\infty} s_i = \lim_{i\to\infty} t_i =\infty \quad
\text{and} \quad (g_1^{-s_i}g_2^{t_i}K) \cap K \neq \emptyset.
\]
In particular, the subgroup of $\Aff(N)$ generated by $g_1$ and $g_2$ does not
act properly discontinuously on $N$.
\end{prop}

\begin{proof}
The last part follows from the group theoretical argument given in
\cite{AMS}, Corollary $2.3$. It is independent of our
generalization.

Choose a norm $\norm{\cdot}$ on $\Ln$. It defines a left-invariant
metric $d$ on $N$.\\ Take a pseudohyperbolic element $g\in
\Omega_0$ and use the notations introduced above. If $x=m\cdot
n^0\cdot n^+ =P_g(x)\cdot n^+\in E_g^+$ and $k\in \N$, then
\begin{align*}
d(g^{-k}x,P_g(g^{-k}x)) & = d(g^{-k}x, P_g(x)\tau(g)^{-k})\\
 & = d(g^{-k}(P_g(x)\cdot n^+), P_g(x)\tau(g)^{-k})\\
 & = d(g^{-k}(P_g(x))\al^{-k}(n^+), P_g(x)\tau(g)^{-k})\\
 & = d( P_g(x)\tau(g)^{-k}\al^{-k}(n^+), P_g(x)\tau(g)^{-k})\\
 & = d(\al^{-k}(n^+), 1)
\end{align*}
since the metric is left-invariant and
$g^{-k}(ab)=g^{-k}(a)\al^{-k}(b)$ for all $a,b\in N$.
Using the
fact that $d(\exp(A),1)\le \norm{A}$ we have
\begin{align*}
d(\al^{-k}(n^+), 1) & \le \norm{\log(\al^{-k}(n^+))}\\
 & = \norm{\al_{\ast}^{-k}(\log(n^+))}\\
 & \le ce^{-bk}\norm{\log(n^+)}
\end{align*}
for some constants $b,c>0$ only depending on $\al_{\ast}$. (Use
that $\alpha_{\ast\mid \Ln^+(g)}^{-1}$ has only eigenvalues
$\lambda$ of modulus $|\lambda|<1$).
Thus we have
\begin{equation}\label{8}
d(g^{-k}x,P_g(g^{-k}x)) \le ce^{-bk}\norm{\log(n^+)}
\end{equation}
for all $k\in \N$ and $x\in E_g^+$.

Fix a point $m(g)$ on $C_g$ and write $C_g=m(g)N^0(g)$. Let
\[
R(g)=\{m(g)\tau(g)^t \mid 0\le t<1 \}.
\]
For every $x\in E_g^+$ there exists a unique integer $k(x,g)$ such
that $P_g(g^{k(x,g)}x)\in R(g)$.

Write $E^+_{g_1}=m(g_1)D^+(g_1)$ and $E^+_{g_2}=m(g_2)D^+(g_2)$.
We claim that $E^+_{g_1}\cap E^+_{g_2}$ is not empty. To show
this, we use the following lemma which can be easily proved by
induction on the nilpotency class of $\Ln$.

\begin{lem}\label{surj}
Suppose that $\Ln$ is a nilpotent Lie algebra which is the sum of
two subalgebras: $\Ln=\La+\Lb$. Let $N=\exp(\Ln)$, $A=\exp(\La)$
and $B=\exp(\Lb)$. Then the map $\phi\colon A\times B\ra N$,
$(a,b)\mapsto ab$ is surjective.
\end{lem}

Since $g_1,g_2$ form a positive pair we have
$\Ln=\Ld^+(g_1)+\Ld^+(g_2)$. Hence by Lemma \ref{surj} we may
write $m(g_1)^{-1}m(g_2)=m_1m_2$ with $m_1\in D^+(g_1)$ and
$m_2\in D^+(g_2)$. Then
\begin{align*}
E^+_{g_2} & = m(g_2)D^+(g_2) = m(g_1)m(g_1)^{-1}m(g_2)D^+(g_2)\\
 & = m(g_1)m_1m_2D^+(g_2) =  m(g_1)m_1D^+(g_2)
\end{align*}
so that $m(g_1)m_1\in E^+_{g_2}$. On the other hand, $m(g_1)m_1\in E^+_{g_1}$.
Hence we have found an element
\[
x_0=m(g_1)m_1\in E^+_{g_1}\cap  E^+_{g_2}.
\]

Now choose a $V\in \Ld^+(g_1)\cap \Ld^+(g_2)$ such that $V\succ_{g_1}0$ and
$V\succ_{g_2}0$. Let $v=\exp(V)$. Fix an element $x_0\in E^+_{g_1}\cap  E^+_{g_2}$
as above.
Then we have for all $i\ge 0$
\[
x_i:=x_0\cdot v^i \in  E^+_{g_1}\cap  E^+_{g_2}.
\]
Let $s_i:=-k(x_i,g_1)$ and $t_i:=-k(x_i,g_2)$. We will make the
computations for $g_1$ and the numbers $s_i$. The argument for
$g_2$ and the numbers $t_i$ is the same.

Let $V_1\in \Ln^0(g_1)$ be so that $V=V_1+W$, where $W\in
\Ln^+(g_1)$. Since $V$ is positive with respect to $g_1$ we have
$V_1=\la_1 T(g_1)$ where $\la_1>0$. Let $v_1=\exp(V_1)$. Note that
$v=v_1\cdot w$ for some $w\in N^+(g_1)$: we can write
\[
\log (v_1\cdot w)=\log(v_1)+\log(w)+W_1
\]
with $W_1\in \Ln^+(g_1)$, since $\Ln^+(g_1)$ is a Lie ideal in
$\Ln^0(g_1)\oplus \Ln^+(g_1)$ and $V=V_1+W$. It follows that
$v^i=v_1^i\cdot w_2$ for some $w_2 \in N^+(g_1)$.

We want to compute $P_{g_1}(g_1^k x_i)$ for any $k\in \Z$. We have
\[
P_{g_1}(g_1^k x_i)=P_{g_1}(g_1^k(x_0\cdot v^i))=P_{g_1}(g_1^k(x_0)\cdot \al_1^k(v^i))
\]
(Of course $\al_1$ denotes the $\Aut(N)$-part of $g_1$). There
exists a $w_3\in N^+(g_1)$ such that
\[
\al_1^k(v^i)=\al_1^k((v_1)^i)\al_1^k(w_2)=v_1^i\cdot w_3=\tau(g_1)^{\la_1 i}w_3
\]
since $V_1=\la_1 T(g_1)$. Writing $x_0=m(g_1)\cdot n^0\cdot n^+$ there exists a
$w_4\in  N^+(g_1)$ such that
\[
g_1^k(x_0)=g_1^k(m(g_1)\cdot n^0)\al_1^k(n^+)=m(g_1)n^0\cdot \tau(g_1)^k\cdot w_4
\]
where $n^0=\tau(g_1)^{r_0}\in N^0(g_1)$.
So we obtain that
\begin{align*}
P_{g_1}(g_1^k x_i) & = P_{g_1}(g_1^k(x_0)\cdot \al_1^k(v^i)) \\
 & = P_{g_1}(m(g_1)n^0\tau(g_1)^k\tau(g_1)^{\la_1i}\cdot \tau(g_1)^{-\la_1i}w_4
\tau(g_1)^{\la_1i}w_3) \\
 & = m(g_1)\tau(g_1)^{r_0+k+\la_1i}.
\end{align*}
This lies in $R(g_1)=\{m(g_1)\tau(g_1)^t \mid 0 \le t <1\}$ if
$0\le r_0+\la_1i+k<1$.
In this case $k$ is the unique integer $k(x_i,g_1)=-s_i$ with this
property. Hence we have $0\le r_0+\la_1i-s_i <1$ for all $i\ge 1$ and
\begin{equation}\label{11}
\lim_{i\to \infty}\frac{i}{s_i}=\lim_{i\to
\infty}\frac{i}{r_0+\la_1i} \le \frac{1}{\la_1}>0.
\end{equation}

Write $x_i=P_{g_1}(x_i)\cdot n_i^+ \in E_{g_1}^+$ with $n_i^+\in N^+(g_1)$.
We have $x_i=x_0\cdot v^i=m(g_1)\cdot n^0\cdot n^+\cdot v^i$. Let
$(W_1,\ldots , W_n)$ be a basis of the nilpotent Lie algebra $\Ln^+(g_1)$.
Using Mal'cev's theorem we can find polynomials $p_1(i), \ldots , p_n(i)$ such that
\[
n^+\cdot v^i=\tau_{g_1}^{\la_1 i}\exp (p_1(i)W_1+\cdots +
p_n(i)W_n).
\]
Hence  $n_i^+ = \exp (p_1(i)W_1+\cdots + p_n(i)W_n)$. So there
exists a polynomial $P(i)$ such that
\[
\norm{\log(n_i^+)}\le P(i).
\]
Using \eqref{8}, \eqref{11} and $b,\la_1 >0$ we obtain
\begin{align*}
0 & \le \lim_{i\to \infty} d(g_1^{-s_i}x_i, P_{g_1}(g_1^{-s_i}x_i)) \le
 \lim_{i\to \infty} c e^{-bs_i}\norm{\log(n_i^+)} \\
  & \le  \lim_{i\to \infty} c e^{-b\la_1 i}P(i)=0.
\end{align*}
It follows that
\[
 \lim_{i\to \infty} d(g_1^{-s_i}x_i, R(g_1))=0.
\]
Hence there exists an upper bound $M_1$ for all these distances. Thus, the
compact set
\[
K_1=\{ x\in N \mid d(x,R(g_1)) \le M_1\}
\]
contains all $g_1^{-s_i}x_i$.

We can obtain in the same way a bound $M_2$ for the distances to
$R(g_2)$, and define the compact set
\[
K=\{ x\in N \mid d(x,R(g_1)) \le M_1\} \cup \{ x\in N \mid d(x,R(g_2))\le M_2\}.
\]
Clearly $g_1^{-s_i}x_i\in K$, $g_2^{-t_i}x_i\in K$ and
$g_1^{-s_i}x_i=(g_1^{-s_i}g_2^{t_i})g_2^{-t_i}x_i$, so that
\[
g_1^{-s_i}g_2^{t_i}K\cap K \neq \emptyset
\]
for all $s_i$ and $t_i$.

\end{proof}

\section{Subgroups of $\Aff(N)$ for $N$ two-step
nilpotent}

In this short section we show that the generalized Auslander conjecture
reduces to the ordinary one if $N$ is two-step nilpotent.
Indeed, if $N$ is two-step nilpotent, a faithful affine
representation
\[
\la\colon \Aff (N)=N\rtimes \Aut (N)\ra \Aff(\R^n)
\]
was constructed in Theorem $4.1$ in \cite{DE2}.  This representation
satisfies the following:
\begin{itemize}
\item Let $i\colon N \hookrightarrow \Aff(N)$ be the embedding given by
$n\mapsto (n,id)$. Then, the  composition $\la\kringel
i\colon N \ra \Aff(\R^n)$ defines a simply transitive action of $N$ on
$\R^n$.  For $n\in N$, $x\in \R^n$ it is given by
\[
n\cdot x  = \la(n,id)(x)\in \R^n.
\]
\item $\la$ maps the subgroup $\Aut(N)$ of $\Aff(N)$ into the subgroup $\GL_n(\R) $
of $\Aff(\R^n)$. It follows that for
every $\alpha\in \Aut(N)$ and for the zero vector $0\in \R^n$, we have that
\[ \la(1,\alpha)(0)=0.\]
\end{itemize}

The following proposition yields the desired reduction of the generalized Auslander
conjecture to the ordinary one:

\begin{prop}\label{3.2}
Let $N$ be a simply connected, connected $2$-step nilpotent Lie group.
Assume that $\Gamma \le \Aff(N)$ acts crystallographically on $N$.
Then $\Gamma$ also admits an affine crystallographic action  on $\R^n$.
\end{prop}

\begin{proof} Let $\la: \Aff(N)\ra\Aff(\R^n)$ be the faithful representation mentioned
above. As $\la$ lets $N$ act simply transitively on $\R^n$, the
evaluation map
\[ e_v:N\ra\R^n:n\mapsto n\cdot0\]
is a diffeomorphism.

Now, $\Aff(N)$ acts on $N$ (via $(n,\alpha)\cdot m=n \alpha(m)$ as
before) and on $\R^n$ (using $\lambda(n,\alpha)$). We can check
that $e_v$ is an $\Aff(N)$--equivariant map, i.e.\ the following
diagram is commutative for any $(n,\alpha)\in \Aff(N)$
$$\xymatrix{
N \ar[r]^{e_v} \ar[d]_{(n,\alpha)\cdot} & \R^n  \ar[d]^{(n,\alpha)\cdot} \\
N \ar[r]_{e_v}  & \R^n  }
$$
Indeed, let $m\in N$, then
\begin{eqnarray*}
(n,\alpha)\cdot e_v(m) &=& \la(n,\alpha) (m\cdot 0)\\
& = & \la(n,\alpha)(\la(m,id)( 0))\\
& = &  \la(n\alpha(m) , \alpha)(0 )\\
& = &  \la(n\alpha(m),id)(\la(1,\alpha)( 0 ))\\
& = & \la(n\alpha(m),id) (0) \\
& = & e_v (n\alpha(m))  \\
& = & e_v ( (n,\alpha)\cdot m)
\end{eqnarray*}
Using this commutative diagram it is now easy to see that a
subgroup $\Gamma$ of $\Aff(N)$ acts crystallographically on $N$,
if and only if it also acts crystallographically (and affinely) on
$\R^n$.
\end{proof}

\section{The conjecture in low dimensions}

In this section we prove that the generalized Auslander conjecture
is true in dimensions $n\le 5$. We have to deal with three cases.
In the first case, the automorphism group of $N$ is already
virtually solvable. Then $\Aff(N)$ is clearly virtually solvable,
and hence all its subgroups $\Gamma$ are virtually solvable. In
the second case, $N$ is $2$-step nilpotent. Then the claim follows
from Proposition $\ref{3.2}$. If neither the first case nor the
second case applies to $N$, the claim is more difficult to prove.
We have to make an appeal to the results of section $\ref{crit}$.
However, in dimension $\leq 5$ there is only one Lie group $N$ for
which this problem arises.

\begin{thm}
Let $N$ be a simply connected and connected nilpotent Lie group of dimension
$n$ with $1\le n\le 5$. Let $\Gamma \le \Aff(N)$ act crystallographically on
$N$. Then $\Gamma$ is virtually polycyclic.
\end{thm}

\begin{proof}
It is enough to show that any such $\Gamma$ is virtually solvable.
Let $\Ln$ be the Lie algebra of $N$. All nilpotent Lie algebras of
dimension $n\le 3$ are nilpotent of class $\leq 2$. If $\dim
\Ln=4$ then $\Ln$ is either of class $\leq 2$ or isomorphic to the
generic filiform Lie algebra $\Ln_4\colon [x_1,x_2]=x_3;
[x_1,x_3]=x_4$. As the derivation algebra of $\Ln_4$ consists of
lower-triangular matrices with respect to this basis, it is
solvable. We can conclude that $\Aut(N_4)$, where
$N_4=\exp(\Ln_4)$, is virtually solvable, which finishes the
argument in dimension 4.

In dimension $5$ we use the list of all nilpotent Lie algebras as
given in \cite{MAG}. It consists of $6$ indecomposable Lie
algebras $\Lg_{5,1}, \ldots , \Lg_{5,6}$ and 3 decomposables. The
algebras $\Lg_{5,1}, \Lg_{5,2}$  and two of the decomposable ones
are nilpotent of class $\leq 2$. The derivation algebra of the
other decomposable one, namely $\Der (\Ln_4\oplus \R$), and the
derivation algebras $\Der (\Lg_{5,3}), \Der (\Lg_{5,5})$ and $\Der
(\Lg_{5,6})$ are clearly solvable. Hence it only remains to
consider the Lie algebra
\[
\Ln = \Lg_{5,4}\; \colon \quad [x_1,x_2]=x_3;\, [x_1,x_3]=x_4; \, [x_2,x_3]=x_5.
\]
This is the free $3$-step nilpotent $2$-generated Lie algebra of dimension $5$.
Let $N=\exp(\Ln)$ and assume that $\Gamma \le \Aff(N)$
acts crystallographically. Assume furthermore that $\Gamma$ is not virtually solvable.
Then also the image of $\Gamma$ inside $\Aut(\Ln)$ under the map
\[
\ell\colon N\rtimes \Aut (N) \ra \Aut(\Ln)
\]
is not virtually solvable. We will show that this leads to a contradiction.\\
A simple calculation shows that, with respect to the basis
$x_1,\ldots, x_5$, $\Aut (\Ln)$ consists of matrices of the form
\[
A_{\al_i,\be_i}=\begin{pmatrix}
\al_1 & \be_1 & 0    & 0        & 0  \\
\al_2 & \be_2 & 0    & 0        & 0   \\
\al_3 & \be_3 & \ga  & 0        & 0   \\
\al_4 & \be_4 & \ast & \al_1\ga & \be_1\ga   \\
\al_5 & \be_5 & \ast & \al_2\ga & \be_2\ga
\end{pmatrix}
\]
where $\ga=\al_1\be_2-\be_1\al_2$ is the determinant of the
$2\times 2$-matrix in the left upper corner (the entries denoted by $\ast$'s are also
determined by the first two columns, but they do not play a role in what follows).
Consider the homomorphism
$\rho\colon \Aut(\Ln)\ra \GL_2(\R)$ given by
\[
A_{\al_i,\be_i} \mapsto
\begin{pmatrix}
\al_1 & \be_1 \\
\al_2 & \be_2
\end{pmatrix}.
\]
It follows that $\rho(\ell(\Gamma))$
is not virtually solvable, because $\ell(\Gamma)$ is not.
Hence, by Tits' alternative $\rho(\ell(\Gamma))$ contains a non-abelian free
subgroup $F_1$. But then, its derived subgroup $F_2=[F_1,F_1]$ is also a
non-abelian free subgroup of $\rho(\ell(\Gamma))$, satisfying
$F_2\subseteq [\GL_2(\R),\GL_2(\R)]=\SL_2(\R)$.
It is well known that for any free non-abelian subgroup $F_2\subseteq  SL_2(\R)$ there
exists an element $g\in F_2$ such that $g$ has no eigenvalues of modulus $1$
(see \cite{CG}). Thus there exists a $g_1\in \Gamma$ such that $\rho (\ell(g_1))\in F_2$
and $\rho(\ell(g_1))$ has no eigenvalue of modulus $1$.
We denote the
eigenvalues of $\rho(\ell(g_1))$ by $\la$ and $1/\la$, with
$\abs{\la}>1$. It is easy to see that there is a basis $(A,B,C,D,E)$ of $\Ln$
with brackets $[A,B]=C,\, [A,C]=D,\, [B,C]=E$ such that
\[
\ell(g_1)=\begin{pmatrix}
\la  & 0     & 0    & 0    & 0  \\
0    & 1/\la & 0    & 0    & 0  \\
0    & 0     & 1    & 0    & 0  \\
\ast & 0     & 0    & \la  & 0  \\
0    & \ast  & 0    & 0    & 1/\la
\end{pmatrix}.
\]
Using the notation of the decompositions \eqref{vsd1} and \eqref{vsd2}
we find that
\begin{align*}
\Ln^0(g_1) & = \langle C\rangle, \quad \Ln^+(g_1)=\langle A,D\rangle  \\
\Ld^+(g_1) & = \langle A,C,D\rangle
\end{align*}
This shows that $T(g_1)=\al C$ for some non-zero $\al$. By rescaling our basis vectors, we may
assume that $T(g_1)=C$.

\medskip

Besides our fixed element $g_1\in \Gamma$, we now choose an
element $h\in \Gamma$ such that $\rho(\ell(h))\in F_2$ and
$\rho(\ell(h))$ does not commute with $\rho(\ell(g_1))$ (this is
possible since $F_2$ is free and non-abelian). We then consider
$g_2=hg_1h^{-1}$. It follows that $\langle \rho(\ell(g_1)),
\rho(\ell(g_2))\rangle $ and hence $\langle \ell(g_1),
\ell(g_2)\rangle$ are free groups. Note that the automorphism
$\ell(g_2)$ has exactly the same eigenvalues as $\ell(g_1)$. Then
there exists a nonzero element of the form $\al A+\be B$ such that
\[
\ell(g_2)(\al A+\be B)=\la (\al A+\be B) \bmod \langle C,D,E \rangle
\mbox{ \ and \ }\ell(g_2)(C) = C \bmod \langle D,E\rangle.
\]
Here we have $\be\neq 0$, otherwise $\langle \ell(g_1),\ell(g_2) \rangle$ would be a
solvable group. Note that
\begin{align*}
\ell(g_2)(\al D + \be E) & = \ell(g_2)([\al A+ \be B, C])\\
& = [\ell(g_2)(\al A+ \be B),\ell(g_2)( C)]\\
 & = [\la (\al A+\be B), C]\\
 & = \la ([\al A+\be B, C])\\
 & =  \la (\al D+\be E).
\end{align*}
It follows that $\dim \Ln^0(g_2)=1$, $\dim \Ln^+(g_2)=2$, so there
exist scalars $\al,\be\neq 0,\ga,\de, \ep, \mu, \nu$ with
\begin{align*}
\Ln^0(g_2) & = \langle C+\mu D+\nu E \rangle  \\
\Ln^+(g_2) & = \langle \al D+ \be E, \al A+\be B+\ga C+\de D+\ep
E\rangle.
\end{align*}
This implies that
\begin{align*}
\Ln & = \Ln^+(g_1)\oplus \Ln^0(g_2) \oplus \Ln^+(g_2) \\
    & = \Ln^+(g_1)\oplus \Ln^0(g_1) \oplus \Ln^+(g_2) .\\
\end{align*}
It follows that $g_1$ and $g_2$ are transversal elements. (Note
that $g_1$ and $g_2$ act fixed-point-free because $\langle g_1,g_2
\rangle$ is a free group and hence torsionfree.) Hence
$\Ld^+(g_1)\cap \Ld^+(g_2)$ is a $1$-dimensional vector space. We
want to show that we can find a positive pair, so that our desired
contradiction follows prom proposition $\ref{prop2.7}$. Let $V\in
\Ld^+(g_1)\cap \Ld^+(g_2)$ be a non-zero vector. This implies that
there are scalars $k,l,m$ and $r,s,t$ such that
\begin{align*}
V & = kA+lC+mD \\
  & = r(\al A+\be B +\ga C+\de D+\ep E) + s(C+\mu D+\ga E)+t(\al D+\be E). \nonumber
\end{align*}
Since $\be\neq 0$ we have $r=k=0$ and $s=l$. For $l=0$ we would obtain $V=0$, hence
we have $l\neq 0$.

\medskip

As $T(g_1)=C$, the semiline $S_{g_1}$ consists of those $V=l C + m D$ with $l>0$.
On the other hand,
$T(g_2)= \xi (C+\mu D + \gamma E)$  for some non-zero $\xi$.
We distinguish two possibilities:
\begin{itemize}
\item If $\xi>0$, then
it is obvious that $S_{g_2}=S_{g_1}$ and hence $g_1$ and $g_2$ form a positive pair.
\item
However, if $\xi$ is negative, we can start all over again and consider
the pair $g_1$ and $g_2^{-1}$. As $T(g_2^{-1})=-T(g_2)$, we obtain that in this case
$S_{g_2^{-1}}=S_{g_1}$ and hence $g_1$ and $g_2^{-1}$ form a postive pair.
\end{itemize}

\end{proof}

\end{document}